\documentclass[letterpaper, 10 pt, conference]{ieeeconf}  

\IEEEoverridecommandlockouts                              

\usepackage{cite}
\usepackage{graphicx}
\usepackage{amsmath,amssymb}
\usepackage{booktabs}
\usepackage{xcolor}
\usepackage[caption=false,font=footnotesize]{subfig}
\definecolor{mattgreen}{RGB}{0, 100, 0} 

\title{\LARGE \bf
Transcription-Induced Failure Modes in 6-DOF Rocket Landing Trajectory Optimization
}

\author{Prayag Sharma$^{1}$, Jonathan Y.M. Goh$^{2}$, Behc\c{e}t A\c{c}\i kme\c{s}e$^{3}$, and Franck Djeumou$^{1}$
\thanks{This work was funded by Toyota Research Institute (TRI)}
\thanks{$^{1}$P. Sharma and F. Djeumou are with MANE Dept., Rensselaer Polytechnic Institute, Troy, NY, USA
        {\tt\small {sharmp6 and djeumf2}@rpi.edu}.}%
\thanks{$^{2}$ J. Goh is with Toyota Research Institute, Los Altos, CA, USA
        {\tt\small jon.goh@tri.global}}%
\thanks{$^{3}$ B. A\c{c}\i kme\c{s}e is with the Department of Aeronautics and Astronautics, University of Washington, Seattle, WA, USA
        {\tt\small behcet@uw.edu}}%
}

\begin{document}

\maketitle
\thispagestyle{empty}
\pagestyle{empty}

\begin{abstract}
Solving optimal control problems via large‑scale NLP solvers depends on discretizing continuous dynamics. Yet, this transcription step hides critical vulnerabilities—most notably truncation error and invariant drift—that can drive solvers toward dynamically infeasible or suboptimal trajectories.
To expose these hidden failures, we introduce a problem- and transcription‑agnostic adversarial objective that leverages the structure of local truncation‑error bounds to aggressively amplify such defects.
When applied to a 6‑DOF rocket‑landing problem, we reveal a stark reliability gap: of fourteen transcription methods tested, only three satisfy rigorous validation criteria. These results also expose a striking performance inversion: even in the absence of classical stiffness, a fourth‑order implicit scheme (GL2) matches the fidelity of a sixth‑order explicit method (RK6). Using B‑series expansions and symplectic Runge–Kutta theorems, we isolate the specific truncation errors and quaternion‑invariant drift responsible for these failures. Crucially, these theoretical vulnerabilities dictate operational performance: in practical lateral‑divert scenarios, the implicit GL2 consistently outperforms the explicit RK6 in both end‑to‑end solve speed and robustness.

\end{abstract}

\section{Introduction}

Optimal control problems (OCPs) are generally solved by transcribing the continuous problem into a finite-dimensional nonlinear programming (NLP) problem\cite{casadi, direct_multiple_shooting}.
While customized successive convexification algorithms like SCvx are tailored for real-time execution \cite{tutorial_SCvx, advances_in_traj_opt}, large-scale NLP solvers such as KNITRO \cite{byrd2006knitro} or IPOPT \cite{wachter2006ipopt} are uniquely suited for rapid offline iteration, dispersion analysis, and rigorous stress testing, as shown in modern benchmarks \cite{Mittelmann_AMPL_NLP_Benchmark_2025}. However, continuous-to-discrete transcription can induce distinct hidden failure modes. When an NLP enforces dynamics via discrete transcription, it assumes negligible local truncation error (LTE) near the optimum. If this assumption is violated, the converged discrete trajectory diverges under continuous open-loop integration. 
Moreover, structural properties of the chosen integrator, e.g. non‑symplecticity, can inadvertently reshape the optimization landscape itself. Thus, reliably exploiting and trusting modern NLP solvers in safety‑critical settings requires that the transcription itself be mathematically rigorous before any solver is applied.
To illustrate how stringent onboard execution limits exacerbate transcription failures, we benchmark these vulnerabilities on a constrained, free-final-time 6-DOF rocket-landing problem, using the same sparse 15-node grid featured in recent NASA SPLICE demonstrations \cite{ElangoTrajoptUtil,ct-scvx,onboard_DQG}. 
Crucially, standard optimal control objectives such as minimum‑fuel landing rarely push the system into regimes that expose these defects, leaving significant reliability gaps untested. To prove and expose the existence of these gaps, we introduce an adversarial objective motivated by the structure of a one-step LTE \cite{LTE_Butcher2008}. By isolating the trajectory-dependent term that dominates standard LTE upper bounds, we construct a principled, transcription-agnostic adversarial objective that amplifies this error. We then derive a tractable approximation of this continuous-time objective and show that max-fuel objective acts as an excellent practical surrogate.


Using this adversarial objective subject to nominal constraints, we evaluate fourteen explicit, implicit, multistep, and symplectic transcriptions with a fixed, robust NLP solver \cite{LTE_Butcher2008,book_non_stiff,book_stiff}. We assess each method against two distinct vulnerabilities: a discrete solution may be dynamically infeasible when rolled out in continuous time, or the transcription’s structure may distort the NLP landscape and drive the solver to a suboptimal minimum. Although all fourteen methods report optimal convergence, only seven satisfy both validation criteria under a minimum-fuel objective. Under the adversarial objective, this number drops to only three (RK6, GL2, and GL3),
revealing substantial hidden reliability gaps and demonstrating the need for adversarial objectives and rigorous validation procedures.

To explain these outcomes, we first trace transcription-induced suboptimality on the adversarial OCP directly to the drift of quadratic invariants. By casting the quaternion dynamics as a quadratic invariant, we show that only the Gauss-Legendre methods used in this study satisfy the symplectic Runge-Kutta coefficient theorems \cite{LTE_Butcher2008} required to prevent such a drift. For non-preserving methods of order below $O(h^6)$, the magnitude of this drift along the adversarial trajectory is large enough to steer the solver towards suboptimal solutions. Although explicitly enforcing quaternion normalization restores convergence for methods with high open-loop accuracy, it fails to rescue standard explicit schemes like RK4 and RK5. 
This analysis identifies invariant drift as a major source of landscape distortion, while also revealing that underlying truncation errors can independently compromise solver reliability.
Motivated by these truncation‑driven failures, we investigate a striking discrepancy in the required integration order: passing our validation checks demands a sixth‑order explicit RK6 method (with or without quaternion normalization), whereas a fourth‑order implicit GL2 method is sufficient. Using a $B$-series expansion over rooted trees, we isolate and characterize the dominant $O(h^5)$ local truncation error (LTE) for the fourth-order RK4 and GL2s schemes. 
After ruling out classical stiffness, further analysis of RK5 shows that the true failure mechanism is not merely per‑step inaccuracy, but the global amplification of these local defects across the transcription. This accumulation makes a sixth‑order explicit method (RK6) necessary to match the fidelity of a fourth‑order implicit scheme.

Crucially, these theoretical vulnerabilities are not confined to contrived, adversarial cases. Rather, they emerge in practical scenarios with safety-critical consequences. To demonstrate this, we evaluate the methods on a standard lateral divert-distance benchmark for rocket landing \cite{G_FOLD,onboard_DQG}. As the vehicle is commanded to increasingly distant touchdown targets, both the optimal flight time and step size naturally grow. This growth forces explicit methods such as RK5 into a computational blow-up, requiring extremely large iteration counts and taking nearly an order of magnitude longer to recover a feasible solution. Although the sixth-order RK6 successfully completes the full sweep, a highly counterintuitive operational reality emerges: it is decisively outperformed by the fourth-order implicit GL2s. By requiring fewer solver iterations and significantly less function-evaluation time, the $3\times$ denser GL2s transcription achieves up to twice the end-to-end solve speed of the explicit RK6.

To the best of the authors' knowledge, this is the first dedicated analysis that systematically exposes transcription-induced failure modes in coarse-mesh 6-DOF rocket-landing NLPs, and rigorously traces the source of these failures back to their principal mathematical mechanisms.

\section{Discrete 6-DOF Rocket Landing Problem}
\label{sec:ocp_discrete}

Since our analysis focuses on the finite-dimensional nonlinear program (NLP) induced by a chosen numerical integration map, we state the rocket-landing problem directly in discrete multiple-shooting form on the normalized horizon \(\tau\in[0,1]\). A detailed continuous-time formulation and background can be found in \cite{ElangoTrajoptUtil,szmuk2017SCvx,ct-scvx,kirk2004optimal}. We partition the normalized horizon into $N-1$ uniform intervals with $N$ nodes $\tau_k=(k-1)h$, $h=1/(N-1)$, $k=1,\dots,N$, and parameterize the control with zero-order hold (ZOH), so that $u(\tau)=u_k$ for $\tau \in [\tau_k,\tau_{k+1})$. For relevance with onboard-flight demonstration literature \cite{onboard_DQG}, as motivated in the introduction, we fix $N=15$ for the remainder of the paper. The discrete decision variables are the state sequence $X:=\{x_k\}_{k=1}^{N},~x_k\in\mathbb{R}^{n_x}$, the control sequence $U:=\{u_k\}_{k=1}^{N-1},~u_k\in\mathbb{R}^{n_u}$, and the time-dilation variable $s>0$, where $t_f=s$ under the normalized-time transformation.



To accommodate both explicit and implicit integration schemes within a unified notation, we introduce stage variables $Y := \{Y_k\}_{k=1}^{N-1}$, where each \(Y_k=[\,y_{k,1}^\top,\dots,y_{k,m}^\top\,]^\top\in\mathbb{R}^{mn_x}\) corresponds to an \(m\)-stage method with stage abscissae \(c_i\in[0,1]\), \(i=1,\dots,m\). For explicit schemes, each \(Y_k\) is determined algebraically by \((x_k,u_k,s)\) and can be eliminated; for implicit schemes, \(Y_k\) is retained as an NLP variable. Let \(\Psi_h(x_k,u_k,s,Y_k)\) denote the one-step numerical map induced by the chosen integrator on \([\tau_k,\tau_{k+1}]\).
For a fixed transcription \(\Psi_h\), let \(\mathcal{A}\) denote the set of all discrete tuples \((X,U,s,Y)\) satisfying the induced discrete dynamics, path constraints, and boundary conditions. The multiple-shooting discrete optimal control problem is then 

\begin{align}\label{eq:ocp_discrete}
\min_{(X,U,s,Y)\in\mathcal{A}} \quad &
J=\Phi_d(x_N,s) + h\sum\nolimits_{k=1}^{N-1} L_d(x_k,u_k,s,Y_k),
\end{align}
where \(\mathcal{A}\) is equivalently characterized by the constraints
\begin{subequations}
\label{eq:general_discrete_ocp_constraints}
\begin{align}
&x_{k+1}=\Psi_h(x_k,u_k,s,Y_k), ~~ k=1,\dots,N-1, \\
&g_d(x_k,u_k,s,Y_k)\le 0, ~ \qquad k=1,\dots,N-1, \\
&b_d(x_0,x_N,s) =0.
\end{align}
\end{subequations}
Here, \(\Phi_d\) and \(L_d\) denote the terminal and stage costs, while \(g_d\) and \(b_d\) collect the discrete path and boundary constraints, respectively. 

\subsection{6-DOF Rocket dynamics and constraints}
The vector field of the system is given by 
\begin{align}
\nonumber &f(x,u,p) =\\
&\begin{bmatrix}
-\alpha_{mdt}\|T_B\|-\beta_{mdt} \\
v_I \\
\frac{1}{m}C_{BI}^\top\!\left(T_B+A_B\right)+g_I \\
\frac{1}{2} q \otimes \begin{bmatrix} \omega_B^\top & 0 \end{bmatrix}^\top \\
J_B^{-1}\!\left(r_{TB}\times T_B+r_{cpB}\times A_B-\omega_B\times(J_B\omega_B)\right)
\end{bmatrix},
\label{eq:rocket_dynamics_discrete_sec}
\end{align}
and the path constraints are imposed through
\begin{equation}
\label{eq:rocket_path_constraints_discrete}
g_d(x_k,u_k,s)=
\begin{bmatrix}
m_{dry}-m_k \\
\|H_\gamma r_{I,k}\| - r_{I,k,x}\tan(\gamma_{gs}) \\
\|H_\theta q_k\| - \sin(\theta_{max}/2) \\
\|\omega_{B,k}\| - \omega_{max} \\
\|T_{B,k}\| - T_{max} \\
\cos(\delta_{max})\|T_{B,k}\| - T_{B,k,x} \\
\|v_{I,k}\| - V_{max} \\
T_{min} - \|T_{B,k}\| \\
s - s_{max} \\
s_{min} - s
\end{bmatrix}
\le 0 .
\end{equation}

Here, \(x=[\,m,\ r_I^\top,\ v_I^\top,\ q^\top,\ \omega_B^\top\,]^\top\in\mathbb{R}^{14}\) is the state, \(m\) is the mass, $r_I$ and $v_I$ are the inertial-frame position and velocity, $q$ is the attitude quaternion, and $\omega_B$ is the body-frame angular velocity. The control is the body-frame thrust vector \(u=T_B\in\mathbb{R}^3\). The aerodynamics force is given by \(A_B:=-\tfrac{1}{2}\rho S_A\|v_I\|(C_A C_{BI})v_I\), where \(C_{BI}=C_{BI}(q)\) is the direction-cosine matrix associated with the quaternion \(q\), and \(\otimes\) denotes quaternion multiplication. 
The plant parameters are \(p=\{\alpha_{mdt},\beta_{mdt},g_I,J_B,r_{TB},r_{cpB},\rho,S_A,C_A\}\). Detailed parameter values and constraint definitions are omitted
for space and can be found in \cite{ElangoTrajoptUtil,szmuk2017SCvx,ct-scvx}.
Under the normalized‑time transformation, the dynamics satisfy 
$$\dot{x}(\tau)=s\,f(x(\tau),u(\tau),p).$$
The nominal control objective is minimum fuel, and is enforced with \(J_d=\Phi_d=-e_m^\top x_N\), where \(e_m=[\,1,0,\dots,0\,]^\top\) extracts the terminal mass component. 


The nodal path constraints are imposed as \(g_d(x_k,u_k,s)\le 0\), \(k=1,\dots,N\), where $H_\gamma=[\,0_{2\times 1}\;\; I_2\,]$, and $H_\theta=[\,0_{2\times 1}\;\; I_2\;\; 0_{2\times 1}\,]$. The components of \(g_d\) enforce, respectively, dry-mass, glide-slope, tilt-angle,
angular-rate, upper-thrust, thrust-pointing, speed, lower-thrust, upper-time, and
lower-time bounds. The boundary conditions enforce the prescribed initial state and terminal touchdown conditions. Specifically, \(x_0\) is fixed to \(x_{\mathrm{init}}=[\,m_{wet},\,r_{I,1}^\top,\,v_{I,1}^\top,\,0,\,0,\,0,\,1,\,0_{1\times 3}\,]^\top\), while the terminal conditions require \(r_{I,N}=r_{I,K}\), \(v_{I,N}=[-0.1,\,0,\,0]^\top\), \(q_N=[0,\,0,\,0,\,1]^\top\), and \(\omega_{B,N}=[0,\,0,\,0]^\top\).

\section{An Adversarial Objective for Stressing Transcription-Induced Errors} \label{sec:worst_case}
This section constructs a modified rocket-landing OCP \eqref{eq:ocp_discrete} that preserves the original boundary and path constraints \eqref{eq:general_discrete_ocp_constraints}, but replaces the nominal fuel-minimization objective with an adversarial cost designed to amplify \emph{one-step} local truncation errors (LTE) on a fixed mesh. The goal is to expose integrator-dependent robustness differences that remain hidden under the nominal objective. The resulting adversarial instance is then used to benchmark multiple integrators under identical constraints,  discretization, and initializations.

\subsection{Flow map, numerical map, and one-step LTE}
\label{subsec:lte_in_nlp}

\noindent\textbf{Definition 1.}
\textit{Let \(\Phi_h(x_k,u_k,s)\) be the exact one-step flow map over \([\tau_k,\tau_{k+1}]\), defined by \(\Phi_h(x_k,u_k,s):=x(\tau_k+h) = x_{k+1}\), where \(x(\cdot)\) solves \(\dot{x}(\tau)=s\,f(x(\tau),u_k,p)\) with initial condition \(x(\tau_k)=x_k\). Let \(\tilde{x}_{k+1}=\Psi_h(x_k,u_k,s,Y_k)\) denote the corresponding one-step numerical map induced by the chosen integrator, so that \(\tilde{x}_{k+1}\) is the discrete state predicted at \(\tau_{k+1}\).}

\medskip
In the direct multiple-shooting transcription, the dynamics over $[\tau_k,\tau_{k+1}]$ are enforced via the defect constraints
\begin{equation}
d_k(x_{k+1},x_k,u_k,s,Y_k):=x_{k+1}-\Psi_h(x_k,u_k,s,Y_k)=0.
\end{equation}
The corresponding \emph{one-step local truncation error} (LTE) is
\begin{equation}
\label{eq:lte_map_exp}
\delta_h(x_k,u_k,s,Y_k):=\Phi_h(x_k,u_k,s)-\Psi_h(x_k,u_k,s,Y_k).
\end{equation}
Thus, even when the NLP satisfies the discrete defect exactly, the propagated state variable \(x_{k+1}\) may still differ from the exact continuous-time flow by \(\delta_h\).

\subsection{LTE upper bounds motivate an adversarial objective}
\label{subsec:lte_bound_objective}

\noindent\textbf{Assumption A.}
\emph{For each interval $[\tau_k,\tau_{k+1}]$ under ZOH control $u_k$, the initial value problem defining $\Phi_h(x_k,u_k,s)$ has a unique solution $x(\cdot)$ that is $(p\!+\!1)$-times continuously differentiable in $\tau$, with $\sup_{\tau\in[\tau_k,\tau_{k+1}]}\|x^{(p+1)}(\tau)\|<\infty$.}


\medskip
\noindent\textbf{Assumption B.}
\emph{The admissible set \(\mathcal{A}\) of decision variables \((X,U,s,Y)\), obtained by enforcing the boundary and path constraints at mesh nodes and, where applicable, at internal stage points, is closed and bounded, and hence compact.}




\medskip
\noindent\textbf{Proposition 1.}
\emph{Let \(\Psi_h\) be an order-\(p\) explicit or implicit Runge--Kutta method and define \(\delta_h=\Phi_h-\Psi_h\) as in \eqref{eq:lte_map_exp}. Let $\|x^{(p+1)}\|_{L^\infty_k}:=\sup_{\tau\in[\tau_k,\tau_{k+1}]}\|x^{(p+1)}(\tau)\|$. Under Assumptions~A--B, there exist method-dependent constants \(\kappa_{\Psi,1},\kappa_{\Psi,2}>0\) such that for every interval \([\tau_k,\tau_{k+1}]\),}
\begin{equation}
\|\delta_h(x_k,u_k,s,Y_k)\|
\le
h^{p+1}\Big(
\kappa_{\Psi,1}\|x^{(p+1)}\|_{L^\infty_k}
+
\kappa_{\Psi,2}\Xi_k
\Big).
\label{eq:lte_bound_sup}
\end{equation}
\emph{where \(\Xi_k \ge 0\) denotes the remaining stage-derivative contributions.}

\medskip
\noindent\textit{Proof.}
For an order-\(p\) Runge--Kutta method, the standard local-error estimate gives
\begin{equation}
\|\delta_h(z_k)\|
\le
h^{p+1}\Big(
\kappa_{\Psi,1}(z_k)\|x^{(p+1)}\|_{L^\infty_k}
+
\kappa_{\Psi,2}(z_k)\Xi_k
\Big).
\end{equation}
with \(z_k:=(x_k,u_k,s,Y_k)\); see \cite[Thm.~II.3.1]{book_non_stiff}. Under Assumption~B, the admissible set is compact, so the method-dependent factors $(\kappa_{\Psi,1}(z_k),\kappa_{\Psi21}(z_k))$ admit finite suprema over \(\mathcal A\). This yields the uniform constants \(\kappa_{\Psi,1},\kappa_{\Psi,2}\) in \eqref{eq:lte_bound_sup}. \hfill\(\square\)

\medskip
\noindent\textbf{Remark 1.}
\textit{ \(\sup_{\tau\in[\tau_k,\tau_{k+1}]}\|x^{(p+1)}(\tau)\|\) in Proposition~1 enters as a positive contributor to the LTE upper bound. We therefore maximize this term as a tractable surrogate for inflating the bound on a fixed mesh.}


\medskip
\noindent\textbf{Remark 2.}
\textit{Proposition~1 shows that the leading LTE upper bound term depends on both the integrator, through \(\kappa_\Psi\) and order \(p\), and the trajectory through \(\sup_{\tau\in[\tau_k,\tau_{k+1}]}\|x^{(p+1)}(\tau)\|\). To decouple integrator and trajectory dependence, we parametrize \(p\) to construct a family of increasingly demanding feasible trajectories, and then examine when different integrators fail.}


\medskip
Motivated by Proposition~1, we consider the ideal adversarial objective for LTE amplification for a fized $p$
\begin{align}
\label{eq:Jsup_def}
\max_{(X,U,s,Y)\in\mathcal{A}} &
J_{\sup}(X,U,s,Y)
:=
\sum_{k=1}^{N-1}
\sup_{\tau\in[\tau_k,\tau_{k+1}]}\|x^{(p+1)}(\tau)\|,
\end{align}
so that maximizing \(J_{\sup}\) over \(\mathcal A\) selects an admissible trajectory that inflates a tractable trajectory-dependent component of the LTE upper bound \eqref{eq:lte_bound_sup}. We next formulate an implementable approximation of this interval wise continuous objective.

\subsection{Supremum approximation via $\ell_r$-smoothing}
\label{subsec:discrete_sup_qnorm}

For any admissible tuple \((X,U,s,Y)\in\mathcal A\), let \(\tau_{k,i}:=\tau_k+c_i h\) denote fixed sampling locations within \([\tau_k,\tau_{k+1}]\), where \(c_i\in[0,1]\), \(i=1,\dots,m\). We then define the nonnegative stage samples $a_{k,i}^{(p)}(x_k,u_k,s):=\|x^{(p+1)}(\tau_{k,i})\|.$, where \(x(\cdot)\) denotes the exact local trajectory on \([\tau_k,\tau_{k+1}]\) generated by the interval data \((x_k,u_k,s)\). Using these stage samples, define the stage-sampled max functional
\begin{equation}
J_{\infty}(X,U,s,Y)
\;:=\;
\sum_{k=1}^{N-1}\max_{i=1,\dots,m} a_{k,i}^{(p)}(x_k,u_k,s),
\label{eq:Jinf_def}
\end{equation}
which satisfies \(J_{\infty}\le J_{\sup}\) on \(\mathcal A\) by construction, since the finite set $\{\tau_{k,i}\}_{i=1}^m \subset [\tau_k, \tau_{k+1}]$.

\medskip
\noindent\textbf{Lemma 1.}
\emph{For \(a_1,\dots,a_m\ge 0\) and \(r\ge 1\),}
\begin{equation}
\max_i a_i \;\le\; \Big(\sum_{i=1}^m a_i^r\Big)^{1/r} \;\le\; m^{1/r}\max_i a_i,
\label{eq:lr_max_bound}
\end{equation}
\emph{since \(m^{1/r}\to 1\) as \(r\to\infty\), \(\big(\sum_{i=1}^m a_i^r\big)^{1/r}\to \max_i a_i\) as \(r\to\infty\).} see \cite[Thm.~1.2.2]{book_applied_math}



\medskip
Following Lemma~1, we define a stage-sampled surrogate of the adversarial objective by replacing the per-interval stagewise maximum with its \(\ell_r\)-approximation,
\begin{equation}
J_r(X,U,s,Y)
\;:=\;
\sum\nolimits_{k=1}^{N-1}
\Big(\sum\nolimits_{i=1}^m a_{k,i}^{(p)}(x_k,u_k,s)^r\Big)^{1/r}.
\label{eq:Jr_def}
\end{equation}

\medskip
\noindent\textbf{Proposition 2.}
\emph{Let \(J_{\sup}\) be defined in \eqref{eq:Jsup_def}, and \(J_{\infty}\) and \(J_r\) in \eqref{eq:Jinf_def}--\eqref{eq:Jr_def}. Then, for every admissible tuple \((X,U,s,Y)\in\mathcal A\),}
\begin{equation}
J_{\infty}\;\le\;J_r\;\le\;m^{1/r}J_{\infty}\;\le\;m^{1/r}J_{\sup}.
\label{eq:prop2_chain}
\end{equation}
\emph{Moreover, \(J_r\to J_{\infty}\) pointwise on \(\mathcal A\) as \(r\to\infty\). Hence, maximizing \(J_r\) over \(\mathcal A\) provides stage-sampled surrogate for enlarging the trajectory-dependent component of the LTE bound in Proposition~1, up to the \(\ell_r\)–\(\ell_\infty\) factor \(m^{1/r}\) and the stage-sampling gap \(J_{\infty}\le J_{\sup}\).}



\medskip
\noindent\textbf{Remark 3.}
\textit{In practice, for a given $r,p$, \(J_r\) is implementable once the stage samples \(a_{k,i}^{(p)}\) are available as functions of the NLP variables. For implicit Runge--Kutta schemes, the required internal stages are naturally lifted into the transcription; for explicit schemes, they are algebraically determined by the node variables. In the present study, GL3 (Table \ref{tab:integrator_properties}) provided sufficient internal resolution.}


\subsection{Fuel maximization as a cheap surrogate for $\ell_r$-objective}
\label{subsec:practical_surrogate_objective}

\begin{figure}[!htbp]
    \centering
    \includegraphics[width=\linewidth]{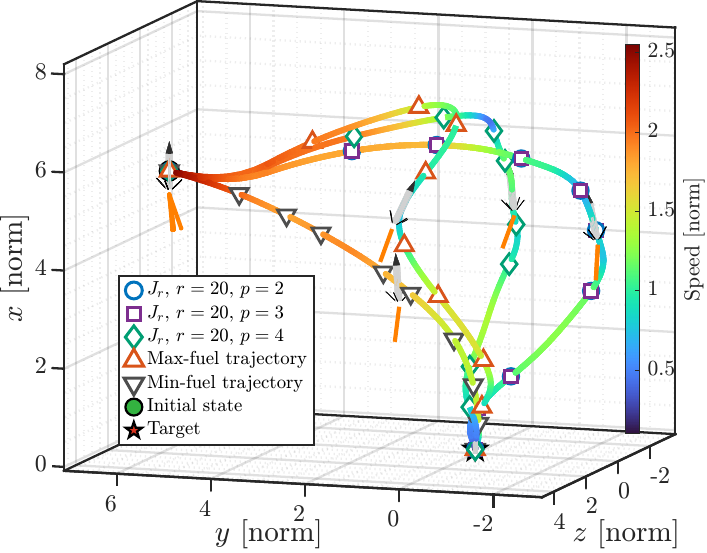}
    \caption{Feasible landing trajectories obtained from $J_r$ for $r=20$ and $p=2,3,4$, and from the proposed \emph{max-fuel} and min-fuel objectives.}
    \label{fig:derivative_max_trajs}
\end{figure}

Although the objective in Section~\ref{subsec:discrete_sup_qnorm} directly targets the trajectory-dependent driver of the LTE bound, evaluating it inside the NLP requires constructing \(x^{(p+1)}(\tau_{k,i})\) and differentiating the resulting composite objective through the full transcription. For the rocket landing problem, this derivative generation can take minutes for a fixed transcription before the NLP is solved, making the approach cumbersome for broad integrator benchmarking (Table \ref{tab:integrator_properties}). To address this, we instead seek a system-specific surrogate guided by domain knowledge of rocket landing. For this problem, maximizing fuel consumption, equivalently minimizing terminal mass,
\begin{equation}
\min\nolimits_{(X,U,s,Y)\in\mathcal A} \quad e_m^\top x_N,
\label{eq:maxfuel_obj}
\end{equation}
drives the trajectory toward sustained high-thrust, highly aggressive feasible maneuvers and, empirically, closely reproduces the behavior induced by the derivative-based adversarial objectives. 

\vspace*{1mm}\noindent\textbf{Comparing max fuel comsumption and $\ell_r$-objective} Fig. \ref{fig:derivative_max_trajs} shows the feasible trajectories generated for different $p$ values for a fixed $r=20$, along with min-fuel and max-fuel trajectories. When evaluated using the metric \(J_r\) for $p=2,3,4$, the max-fuel trajectory attains approximately \(91\%,90\%,78\%\) of the optimal objective value, whereas the nominal min-fuel trajectory attained only \(5\%,1\%,0.3\%\) respectively. Thus, although this surrogate is problem-dependent and need not generalize to other systems, it provides a computationally cheap and highly effective adversarial benchmark for the present rocket-landing study.

\section{Integrators in-loop Optimization Benchmark}
Table~\ref{tab:integrator_properties} lists the \(14\) integration methods used in this study. The set is intentionally heterogeneous, to isolate which integrator attributes most strongly influence NLP robustness. Standard background on stiffness, stability, and integrator classes can be found in \cite{book_stiff,book_non_stiff}.

\begin{table*}[!t]
\renewcommand{\arraystretch}{1.1}
\setlength{\tabcolsep}{4pt}
\caption{Properties, computational footprint, and benchmark outcomes of the integration methods. Vars and NNZ are computed for the \(N=15\) 6-DOF rocket-landing transcription in \eqref{eq:ocp_discrete}.}
\label{tab:integrator_properties}
\centering
\footnotesize
\begin{tabular}{c l l c c c c c c|c c|c c|c}
\toprule
\multicolumn{9}{c|}{} & \multicolumn{2}{c|}{\textbf{Min-Fuel}} & \multicolumn{2}{c|}{\textbf{Max-Fuel}} & \multicolumn{1}{c}{\textbf{Max-Fuel}} \\
\cmidrule(lr){10-11} \cmidrule(lr){12-13} \cmidrule(lr){14-14}
\textbf{\#} & \textbf{Method} & \textbf{Type} & \textbf{Order} & \textbf{Stages} & \textbf{Symplectic} & \textbf{A-stable} & \textbf{Vars} & \textbf{NNZ} & \textbf{Optim.} & \textbf{OL} & \textbf{Optim.} & \textbf{OL} & \textbf{OL} \\
\midrule
\textbf{1}  & BDF4                     & Implicit & 4 & 1 & $\times$         & $\times^{\dagger}$ & 253 & 4,091  & $\times$     & $\times$     & $\times$     & $\times$     & $\times$     \\
\textbf{2}  & BDF6                     & Implicit & 6 & 1 & $\times$         & $\times^{\dagger}$ & 253 & 4,425  & $\times$     & $\times$     & $\times$     & $\times$     & $\times$     \\
\textbf{3}  & Implicit Trapezoidal     & Implicit & 2 & 1 & $\times$         & $\checkmark$       & 253 & 4,542  & $\times$     & $\times$     & $\times$     & $\times$     & $\times$     \\
\textbf{4}  & RK38                     & Explicit & 4 & 4 & $\times$         & $\times$           & 253 & 5,041  & $\checkmark$ & $\checkmark$ & $\times$     & $\times$     & $\times$     \\
\textbf{5}  & RK4                      & Explicit & 4 & 4 & $\times$         & $\times$           & 253 & 5,041  & $\checkmark$ & $\checkmark$ & $\times$     & $\checkmark$     & $\times$     \\
\textbf{6}  & RK5 (DoPri5)             & Explicit & 5 & 6 & $\times$         & $\times$           & 253 & 5,041  & $\checkmark$ & $\checkmark$ & $\times$     & $\times$     & $\times$     \\
\textbf{7}  & RK6 (Luther)             & Explicit & 6 & 7 & $\times$         & $\times$           & 253 & 5,041  & $\textcolor{blue}{\checkmark}$ & $\textcolor{blue}{\checkmark}$ & $\textcolor{blue}{\checkmark}$ & $\textcolor{blue}{\checkmark}$ & $\textcolor{blue}{\checkmark}$ \\
\textbf{8}  & AVF 2s                   & Implicit & 2 & 2 & $\checkmark^{*}$ & $\checkmark$       & 253 & 5,858  & $\times$     & $\times$     & $\times$     & $\times$     & $\times$     \\
\textbf{9}  & AVF 3s                   & Implicit & 2 & 3 & $\checkmark^{*}$ & $\checkmark$       & 253 & 5,858  & $\times$     & $\times$     & $\times$     & $\times$     & $\times$     \\
\textbf{10} & Implicit Midpoint (GL1)  & Implicit & 2 & 1 & $\checkmark$     & $\checkmark$       & 253 & 5,858  & $\times$     & $\times$     & $\times$     & $\times$     & $\times$     \\
\textbf{11} & TR-BDF2                  & Implicit & 2 & 2 & $\times$         & $\checkmark$       & 449 & 7,552  & $\times$     & $\times$     & $\times$     & $\times$     & $\times$     \\
\textbf{12} & Gauss--Legendre 2s (GL2) & Implicit & 4 & 2 & $\checkmark$     & $\checkmark$       & 645 & 15,526 & $\textcolor{blue}{\checkmark}$ & $\textcolor{blue}{\checkmark}$ & $\textcolor{blue}{\checkmark}$ & $\textcolor{blue}{\checkmark}$     & $\textcolor{blue}{\checkmark}$ \\
\textbf{13} & Lobatto IIIA 3-stage     & Implicit & 4 & 3 & $\times$         & $\checkmark$       & 841 & 24,318 & $\checkmark$ & $\checkmark$ & $\times$     & $\checkmark$ & $\checkmark$ \\
\textbf{14} & Gauss--Legendre 3s (GL3) & Implicit & 6 & 3 & $\checkmark$     & $\checkmark$       & 841 & 27,258 & $\textcolor{blue}{\checkmark}$ & $\textcolor{blue}{\checkmark}$ & $\textcolor{blue}{\checkmark}$ & $\textcolor{blue}{\checkmark}$ & $\textcolor{blue}{\checkmark}$ \\
\bottomrule
\multicolumn{14}{l}{\footnotesize $^{*}$ Structure-preserving ~~ $^{\dagger}$ \(A(\alpha)\)-stable. ~~\textbf{NNZ}: Total Jacobian plus Hessian non-zeros.~~ \textbf{Optim.}: Converged to the best local minimizer found.}\\
\multicolumn{14}{l}{\footnotesize \textbf{OL}: Open-loop replay satisfies \(\epsilon_{\mathrm{OL}} \le 10^{-2}\).~~ \textbf{Ref. OL}: A numerical map's replay of the best max fuel local minimizer satisfying \(\epsilon_{\mathrm{OL}} \le 10^{-2}\).} \\
\end{tabular}
\end{table*}

\vspace*{1mm}\noindent\textbf{Optimizer choice:} Optimizer selection matters here since the solver must be robust enough to support benchmarking integrator‑induced effects. Since solver performance depends on initialization, scaling, and tuning, a comprehensive optimizer study is beyond the scope of this work. We therefore report only the considerations that motivated our solver choice for the benchmarks.
In recent Mittelmann's AMPL--NLP benchmark \cite{Mittelmann_AMPL_NLP_Benchmark_2025} (47 instances, 2\,hr time limit, benchmark date: 20 Dec 2025), \texttt{KNITRO} compares favorably with common solvers like \texttt{IPOPT} and \texttt{SNOPT}, achieving a smaller shifted geometric mean runtime and solving more instances (47/47 for \texttt{KNITRO}, versus 46/47 for \texttt{IPOPT} and 30/47 for \texttt{SNOPT}). We observed similar behavior on the \(N=15\) rocket-landing problem \eqref{eq:ocp_discrete} subject to \eqref{eq:general_discrete_ocp_constraints}: across multiple transcriptions, initializations, both min- and max-fuel objectives, \texttt{KNITRO} was markedly more reliable than \texttt{IPOPT} and \texttt{SNOPT} in recovering both feasible and optimal solutions. Thus, we fix Interior Point (IP) method based \texttt{KNITRO} as the preferred NLP solver for the integrator study.

\vspace*{1mm}\noindent\textbf{Transcription computational cost:} Under a fixed interior-point solver, per-iteration cost is dominated by nonlinear function/derivative evaluation and sparse linear algebra in the IP step. Both are shaped by the transcription induced by the integrator through the defect constraints and their derivatives. Table~\ref{tab:integrator_properties} reports each integrator's computational footprint for the fixed \(N=15\) rocket landing problem: (a) the total number of NLP decision variables, and (b) the aggregate sparsity (total Jacobian plus Hessian nonzeros). This also captures another tradeoff studied here: implicit stage-based schemes typically add variables and nonzeros through lifted stage constraints, whereas explicit schemes avoid stage lifting but may require deeper computational graphs. Hence per-iteration cost and total iteration counts/runtimes cannot be inferred apriori and must be assessed empirically (see further analysis in Section \ref{subsec:divert_benchmark}).

\vspace*{1mm}\noindent\textbf{Evaluation methodology:} We benchmark 14 integrators on the nominal min-fuel and stress-test max-fuel problems. A best local minimizer is identified using high-fidelity transcriptions (Gauss--Legendre 3s and a 12-stage explicit Runge--Kutta method), which converge to the same solution. Each integrator is evaluated on two criteria: convergence to this minimizer (Optim.) and the open-loop (OL) fidelity of its optimized trajectory.  For the OL test, the optimized controls are propagated through the dynamics \eqref{eq:rocket_dynamics_discrete_sec} via MATLAB's adaptive-step \texttt{ode113} to compute the terminal error $\epsilon_{\mathrm{OL}}=\|r_I^{\mathrm{OL}}(t_f)-r_I^{\mathrm{NLP}}(t_f)\|_2$, requiring $\epsilon_{\mathrm{OL}}\le 10^{-2}$ in normalized coordinates. 

\vspace*{1mm}\noindent\textbf{$\mathbf{3/14}$ successful transcriptions on the adversarial test:} For the nominal minimum-fuel problem, Table~\ref{tab:integrator_properties} shows \(7/14\) methods satisfy both these tests. Because nearly all methods of fourth order or higher succeed (excluding multi-step BDF-4,6), we conclude that a \textit{fourth-order one-step} transcription appears sufficient for the min-fuel problem. Yet, the max-fuel stress test exposes this baseline's inadequacy, leaving only RK6, GL2s, and GL3s as successful transcriptions.

\vspace*{1mm}\noindent\textbf{Multiple mechanisms distorting the optimization landscape:} To isolate the impact of a transcription's inherent numerical inaccuracy on the optimization outcome, we introduce the Ref. OL test (Table~\ref{tab:integrator_properties}). This test assesses a map's ability to faithfully propagate the best-known optimal controls. Analyzing the failures of RK4 and Lobatto~IIIA provides critical insight into how transcriptions distort the optimization space. RK4 fails to propagate the best-known optimal controls despite its own trajectory passing replay, suggesting its coarse-step map struggles with the aggressive reference trajectory and distorts the feasible landscape toward artificially benign solutions. In contrast, Lobatto~IIIA accurately propagates both its own trajectory and the optimal controls, yet still steers the NLP to a different local minimizer. This complementary failure proves a critical point: open-loop propagation accuracy does not uniquely dictate optimization landscape distortions, implying the presence of supplementary factors which we analyze in Section \ref{subsec:quat_norm_invariance}.

\vspace*{1mm}\noindent\textbf{Integrator's order alone is an unreliable predictor of numerical accuracy:} The results in Table~\ref{tab:integrator_properties} demonstrate that nominal order alone is an unreliable predictor of accuracy. While RK4 and GL2 are both fourth-order, RK4 suffers from significantly larger open-loop replay errors, necessitating a sixth-order explicit method to achieve parity with GL2. A systematic analysis of the stiffness ratio across various min- and max-fuel trajectories yields values strictly below 100, placing the dynamics well outside the traditionally stiff regime $\ge O(10^3)$ \cite{book_stiff}. Consequently, the failures observed under maximum stress cannot be attributed to numerical stiffness. We analyze the underlying mechanisms driving these outcomes in Section \ref{subsec:lte_gl2_rk4}.

\section{Principled analysis}
\label{sec:principled_analysis}


The following subsections deconstruct the mechanisms behind failure modes observed in our stress tests.
We explain how a transcription may distort the optimization landscape irrespective of its open-loop accuracy, and demonstrate how the magnitude and accumulation of local truncation error (LTE) can dictate the selection of a preferred method.

\subsection{Transcription-induced suboptimality via quaternion drift}
\label{subsec:quat_norm_invariance}

The primary mechanism found is a \emph{drift in the quadratic invariant}. In continuous time, the quaternion dynamics preserve the unit-norm constraint $q^\top q = 1$ exactly. However, non-sympletic discrete integration schemes do not.
When the NLP solver enforces feasibility at tolerances much tighter than the integrator’s natural drift scale, the optimization landscape is effectively altered. This leads to termination at stationary points that meet the optimization feasibility tolerance, hence significantly narrowing the solution space.

Let $q(t)\in\mathbb{R}^4$ denote the scalar-last quaternion, and let $\omega(t)\in\mathbb{R}^3$ be the body angular rate. The quaternion kinematics in \eqref{eq:rocket_dynamics_discrete_sec} can be written as
\begin{equation}
\dot q(t) = 0.5 \cdot \Omega(\omega(t))\,q(t)~~\text{with}~~ \Omega(\omega)^\top = -\Omega(\omega),
\label{eq:q_kinematics_matrix}
\end{equation}
where $\Omega(\omega)\in\mathbb{R}^{4\times 4}$ is a matrix representation of quaternion multiplication by $[\omega^\top,\,0]^\top$ under the adopted convention. A standard property of $\Omega(\omega)$ is skew-symmetry:

\medskip
\noindent\textbf{Lemma 2.}
\emph{If $q(\cdot)$ satisfies \eqref{eq:q_kinematics_matrix}, then $q^\top q$ is invariant. In particular, if $q(t_0)^\top q(t_0)=1$, then $q(t)^\top q(t)=1$ for all $t$.}

\noindent\emph{Proof.}
Differentiate $q^\top q$ and use \eqref{eq:q_kinematics_matrix}:
\[
\frac{d}{dt}(q^\top q)=2 q^\top \dot q
= q^\top\big(\Omega(\omega)+\Omega(\omega)^\top\big)q
=0 \eqno{\square}
\]


\medskip
\noindent\textbf{Definition 2.}
\textit{Consider $\dot x = f(x)$ and a quadratic functional $I(x) = x^\top M x$ where $M=M^\top$. $I(\cdot)$
is invariant in continuous time if $\nabla I(x)^\top f(x)=0$ for all $x$ in the domain of interest.}

\medskip
\noindent\textbf{Theorem 2.}\label{thm_2}
\emph{Let a Runge--Kutta method with Butcher tableau $(A,b)$ be applied to $\dot x=f(x)$. If the coefficients satisfy the symplecticity conditions}
\begin{equation}
b_i a_{ij} + b_j a_{ji} - b_i b_j = 0
\qquad \forall\, i,j,
\label{eq:sympl_condition}
\end{equation}
\emph{then the method preserves every quadratic invariant $I(x)=x^\top M x$ exactly: $I(x_{n+1})=I(x_n),~ \forall x : \nabla I(x)^\top f(x)=0$}

\noindent\emph{Proof .}
This is a classical result for symplectic Runge--Kutta schemes (Refer to \cite[Thm.~4.3]{book_symplectic}) the condition \eqref{eq:sympl_condition} implies discrete conservation of all quadratic invariants. \hfill$\square$

\medskip
For the rocket model \eqref{eq:rocket_dynamics_discrete_sec}, Theorem~2 applies to the quaternion subsystem \eqref{eq:q_kinematics_matrix} since Lemma~2 establishes that $I(q)=q^\top q$ is an invariant. Consequently, symplectic Runge--Kutta schemes preserve $q^\top q$ \emph{exactly up to roundoff}. We further compute this symplecticity residual $S_{ij}(A,b) := b_i a_{ij} + b_j a_{ji} - b_i b_j$, $(i,j=1,\dots,s)$ for the integrators under evaluation. A method is symplectic if and only if $S_{ij}\equiv 0$.


\medskip
\noindent\textbf{(i) Gauss--Legendre GL2 (2-stage, order 4) and GL3 (3-stage, order 6).}
Using the coefficients of GL2
\[
A=\begin{bmatrix}
\frac14 & \frac{3-2\sqrt3}{12}\\[2pt]
\frac{3+2\sqrt3}{12} & \frac14
\end{bmatrix},
\qquad
b=\begin{bmatrix}\frac12\\[2pt]\frac12\end{bmatrix},
\]
we have by direct substitution that $S_{ij}(A,b)=0$ for all $i,j$. One similarly obtains $S_{ij}(A,b)=0$ for $(A,b)$ corresponding to GL3 \cite{book_stiff}. Hence GL-2,3 are symplectic and preserve quadratic invariants such as $q^\top q$.

\medskip
\noindent\textbf{(ii) Lobatto IIIA (3-stage, order 4).}
Using the coefficients
\[
c=\begin{bmatrix}0\\[2pt]\frac12\\[2pt]1\end{bmatrix},
\qquad
b=\begin{bmatrix}\frac16\\[2pt]\frac23\\[2pt]\frac16\end{bmatrix},
\qquad
A=\begin{bmatrix}
0 & 0 & 0\\[2pt]
\frac{5}{24} & \frac13 & -\frac{1}{24}\\[2pt]
\frac16 & \frac23 & \frac16
\end{bmatrix},
\]
the symplecticity residual evaluates to $S(A,b)\neq 0$. Therefore Lobatto IIIA is \emph{not} symplectic as a standalone method (in contrast to certain partitioned Lobatto pairs \cite{book_stiff}) and it does not inherit the quadratic-invariant preservation guarantee of Theorem~1. \textit{This predicts quaternion norm drift regardless of its implicit and $A$-stable nature (Table \ref{tab:integrator_properties}).}

\medskip
\noindent\textbf{(iii) Explicit Runge--Kutta methods (general case).}
An $s$-stage explicit Runge--Kutta (RK) method has a strictly lower-triangular Butcher matrix $A$, hence $a_{ii}=0$ for all $i$.

\medskip
\noindent\textbf{Proposition 3.}
\emph{Let $(A,b)$ define the coefficients of a consistent explicit Runge--Kutta method (i.e., $A$ be lower triangular and $\sum_{i=1}^s b_i = 1$). Then, $(A,b)$ do not satisfy the symplecticity conditions \eqref{eq:sympl_condition}}.

\noindent\emph{Proof.}
Set $i=j$ in \eqref{eq:sympl_condition} to obtain the necessary condition
\[
b_i a_{ii} + b_i a_{ii} - b_i^2 = 0
\quad\Longrightarrow\quad
2b_i a_{ii} = b_i^2.
\]
For an explicit RK method, $a_{ii}=0$, hence $b_i^2=0$ and therefore $b_i=0$ $\forall i$. This contradicts consistency $\sum_i b_i=1$. Thus no consistent explicit RK method can satisfy \eqref{eq:sympl_condition}. \hfill$\square$

\medskip
\noindent\textbf{Corollary 1.}
\emph{All explicit RK schemes used in this study (including RK38, RK4, RK5, RK6) are not symplectic, and therefore do not guarantee the exact quadratic-invariant preservation of Theorem~2. In particular, they should generically exhibit step-size and truncation error dependent drift in the quaternion norm $q^\top q$ (beyond roundoff).}

\begin{figure}[!bhtp]
\centering
\includegraphics[width=1\linewidth]{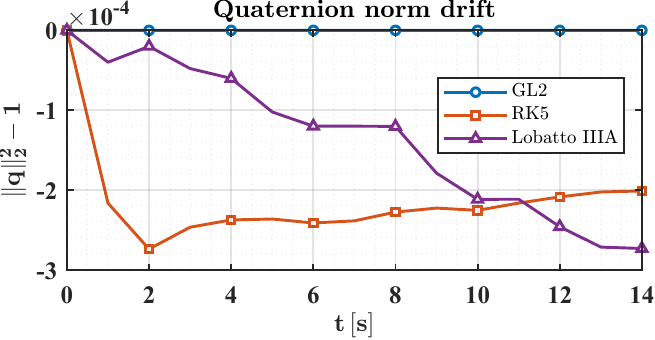}
\caption{Quaternion norm drift $\|q(t)\|_2-1$ under open-loop propagation on the adversarial trajectory. GL2 exhibits near-zero drift around $10^{-12}$, while RK5 and Lobatto~IIIA drift on the order of $10^{-4}$.}
\label{fig:qnorm_drift}
\end{figure}

\medskip
\vspace*{1mm}\noindent\textbf{Quantifying quaternion norm drift along stress trajectory:} Figure~\ref{fig:qnorm_drift} reports the measured drift $\|q(t)\|_2-1$ under open-loop replay using different transcriptions on the best-known max-fuel trajectory. Consistent with Theorem~2 and the symplecticity checks, GL2 exhibits near-zero drift (up to roundoff $10^{-12}$), whereas non-symplectic maps such as RK5 and Lobatto~IIIA show drift on the order of $10^{-4}$.

\begin{figure}[!tbp]
\centering
\includegraphics[width=1\linewidth]{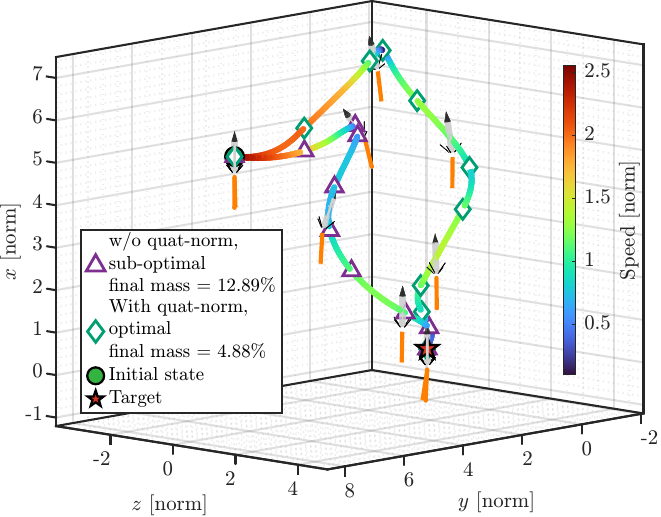}
\caption{Quaternion norm drift impact on Lobatto 3a optimization outcome for the minimum terminal mass (max-fuel) stress test. Quat. normalization producing optimal and $2.6\times$ lower objective value.}
\label{fig:qnorm_lob_3a_3d}
\end{figure}

\vspace*{1mm}\noindent\textbf{Quaternion normalization and its failure modes:} To isolate the role of this invariant drift on the optimizer outcome, we consider the \emph{projected} one-step map obtained by composing the integrator update with quaternion renormalization on the quaternion block, i.e.,
$\bar{\Psi}_h := P\circ \Psi_h$ with $P(q)=q/\|q\|$.
Using $\bar{\Psi}_h$ in the multiple-shooting defect constraints yields convergence for Lobatto~IIIA to the best known local optimum as shown in Fig. \ref{fig:qnorm_lob_3a_3d}. Both optimization outcomes look significantly different with quaternion normalization producing a $2.6\times$ reduction in the objective value. This confirms the interpretation that the observed Lobatto~IIIA suboptimality was primarily driven by the quaternion norm drift.
\emph{However, this projection does not remedy the behavior of standard explicit schemes such as RK4 and RK5. In these methods, we show that the truncation errors themselves are large enough to disrupt the optimization process, independent of the invariant drift.}

\subsection{Characterizing the effect of LTE: GL2 vs. explicit RK4}
\label{subsec:lte_gl2_rk4}


We first characterize the local truncation error (LTE) of RK4 and GL2s along the best-known max-fuel trajectory. Because both are fourth-order methods, their leading LTE is $O(h^5)$. Using a $B$-series expansion over rooted trees \cite[Ch.~3]{butcher_ode_book}, the one-step error in a common order-5 elementary-differential basis $F(t)$ is expressed as:
\begin{equation}
\delta_h(x)
=
h^5 \sum_{r(t)=5}
\frac{1}{\sigma(t)}
\left(
\Phi(t)-\frac{1}{\gamma(t)}
\right)
F(t)(x)
+ O(h^6),
\label{eq:lte_bseries_order5}
\end{equation}
where $r(t)$ is the tree order, $\sigma(t)$ its symmetry factor, $\gamma(t)$ the exact-flow density, and $\Phi(t)$ the method-dependent elementary weight. Thus, under this notation, RK4 and GL2 differ strictly through their tableau-induced weights $\Phi(t)$.

\begin{figure}[!tbp]
    \centering
    \includegraphics[width=\columnwidth]{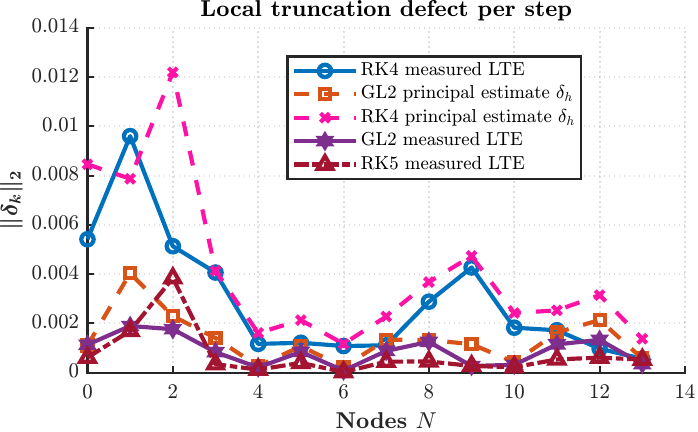}
    \caption{Isolated per-step local truncation errors and principal $O(h^5)$ estimates. While RK4 exhibits significantly larger local defects than GL2, RK5 maintains a comparable magnitude—demonstrating that one-step accuracy alone does not guarantee open-loop fidelity.}
    \label{fig:lte_step_compare}
\end{figure}


To isolate these one-step LTE's for RK4 and GL2, we reset the ground-truth state after every step $h$, allowing us to compute both the analytical $O(h^5)$ estimates in \eqref{eq:lte_bseries_order5} and the true per-step LTE. As shown in Fig.~\ref{fig:lte_step_compare}, the normalized full-state LTE ($\|\delta_k\|_2$) of RK4 is significantly larger than that of GL2, approaching a factor of five at early nodes. The gap between the $h^5$ estimate and the measured LTE highlights the non-negligible impact of higher-order remainders.

\vspace*{1mm}\noindent\textbf{One-step LTE accuracy is not enough to conclude:} However, similar calculations for RK5 (Fig.~\ref{fig:lte_step_compare}) yield per-step LTE magnitudes comparable to GL2, yet RK5 still exhibits severe terminal position errors under continuous-time trajectory replay. Hence, one-step accuracy does not uniquely determine open-loop accuracy. Further analysis reveals that local defects are \textit{disproportionately amplified} by each numerical map over $N$ nodes. Minimizing this global error (GTE) amplification ultimately necessitated the use of the sixth-order RK6. A formal analysis of these global error-dynamics is reserved for future work.

\begin{figure}[htbp]
    \centering
    \includegraphics[width=\linewidth]{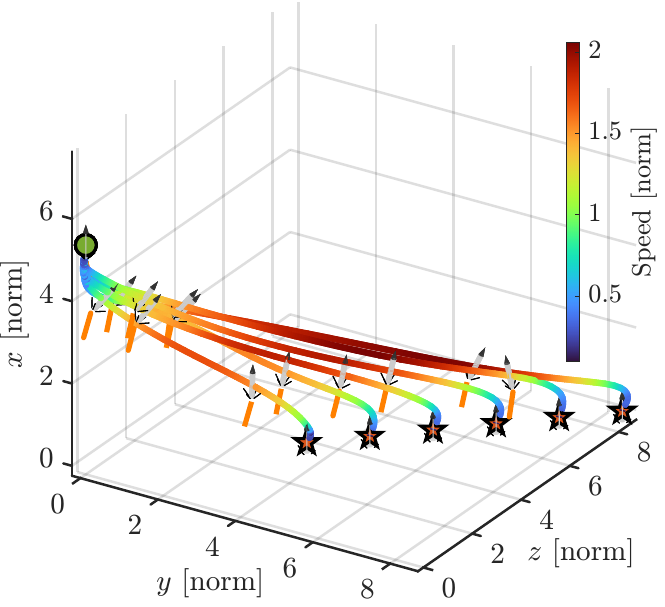}
    \caption{Representative divert maneuvers (normalized coordinates) with different divert landing locations.}
    \label{fig:divert_benchmark}
\end{figure}

\subsection{Operational implications: Lateral divert distance feasibility benchmark}
\label{subsec:divert_benchmark}




    


    

\begin{figure}[htbp]
    \centering
    \includegraphics[width=\linewidth]{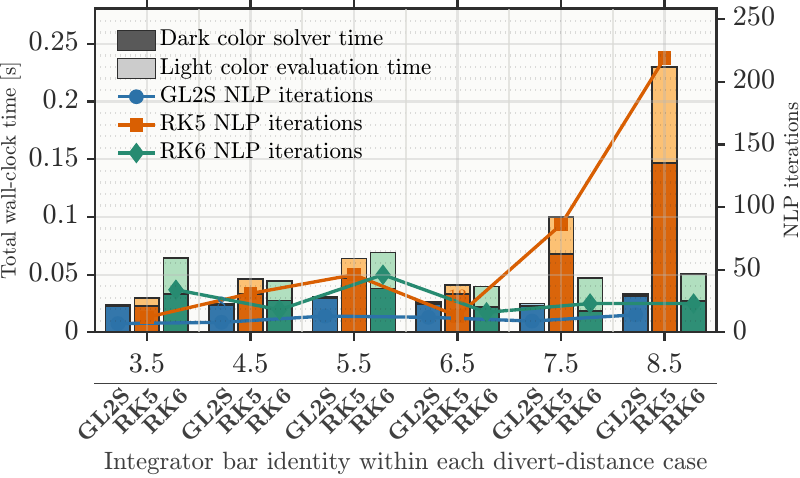}
    \caption{Average wall-clock time (function evaluation + remaining solver time) and iterations across six divert maneuvers.}
    \label{fig:divert_cost_bar}
\end{figure}

\noindent

While the preceding analysis illustrates how unchecked transcriptions can distort offline assessments of vehicle capabilities, these vulnerabilities pose an equally serious threat to large-scale simulation times and real-time trajectory optimization performance. To demonstrate these operational consequences, we evaluate a standard lateral divert benchmark \cite{G_FOLD,onboard_DQG} for a vertically landing rocket across six landing locations (Fig. \ref{fig:divert_benchmark}). Because real-time trajectory recovery supersedes strict optimality during such emergencies, we formulate this maneuver as a \emph{feasibility NLP}, i.e., zero objective function subject to boundary and path constraints.

As the divert distance increases (Fig. \ref{fig:divert_benchmark}), the longer flight times enlarge the physical step size $\Delta t = s \cdot h$ on the fixed $N=15$ mesh.
This amplification makes the integrator increasingly sensitive to truncation and invariant‑drift effects.
This resulting growth forces explicit methods such as RK5 into a computational blow-up (Fig.~\ref{fig:divert_cost_bar}): at larger diverts, it requires prohibitive iteration counts and takes nearly an order of magnitude longer than GL2 to recover a feasible solution.

In contrast, while GL2 and RK6 both remain feasible across all six cases (Fig.~\ref{fig:divert_cost_bar}), a highly counterintuitive reality emerges: Despite generating a substantially larger ($\times3$) derivative footprint (Table~\ref{tab:integrator_properties}), GL2 requires significantly less function-evaluation time and fewer solver iterations. Consequently, the denser GL2 transcription achieves up to twice the end-to-end solve speed of RK6. Moreover, GL2 delivers this performance while maintaining similar iteration counts across maneuvers. Hence, GL2 emerges as the most suitable transcription choice across all presented studies.

\section{Conclusion}


This paper demonstrates that unchecked transcriptions can undermine both trajectory fidelity and optimization performance in coarse-mesh 6-DOF rocket-landing NLPs. The proposed adversarial objective exposed hidden failures that nominal tests did not reveal, reducing fourteen candidate methods to only three successful ones under rigorous validation metrics. A principled analysis further revealed that these failures arise from both invariant drift in non-symplectic methods and the amplification of truncation defects across certain discrete transcriptions. Most importantly, we show that the fourth-order Gauss--Legendre method is more reliable and computationally efficient than similar and higher-order explicit alternatives, such as RK-4,5,6, making it the strongest overall transcription in our benchmarks. Future work will focus on a detailed characterization of global error-dynamics, the impact of NLP solver choice, and broader statistical validation across larger divert scenarios.

\addtolength{\textheight}{-12cm}   


\bibliographystyle{IEEEtran}
\bibliography{bib}

\end{document}